%% file: main.tex
\documentclass{article}

\usepackage[english]{babel}

\usepackage[letterpaper,top=2cm,bottom=2cm,left=3cm,right=3cm,marginparwidth=1.75cm]{geometry}

\usepackage{amsmath,amssymb,amsthm}
\usepackage{graphicx}
\usepackage{bm}
\usepackage[colorlinks=true, allcolors=blue]{hyperref}
\numberwithin{equation}{section}

\newtheorem{theorem}{\bf{Theorem}}[section]
\newtheorem{lemma}{\bf {Lemma}}[section]
\newtheorem{define}{\bf{Definition}}[section]

\newtheorem{proposition}{\bf{Proposition}}[section]
\newtheorem{remark}{\bf{Remark}}[section]

\newcommand{\mbE}{\widehat{\mathbb{E}}}
\newcommand{\mbe}{\widehat{\mathcal{E}}}
\newcommand{\V}{\mathbb{V}}
\newcommand{\mv}{\mathcal{V}}
\newcommand{\sles}{(\Omega,\mathcal{H},\mbE)}

\newcommand{\Vstar}{\widehat{\mathbb{V}}^*}
\newcommand{\mvstar}{\widehat{\mathcal{V}}^*}
\newcommand{\hV}{\widehat{\mathbb{V}}}
\newcommand{\mhV}{\widehat{\mathcal{V}}}
\newcommand{\bE}{\breve{\mathbb{E}}}
\newcommand{\be}{\breve{\mathcal{E}}}
\newcommand{\usigma}{\bar{\sigma}}
\newcommand{\lsigma}{\underline{\sigma}}
\def\tildeCapc{\mathbb C^{\ast}}

\title{On the Law of the Iterated Logarithm for $m$-dependent stationary random variables under sub-linear expectations}
\author{Wang-Yun Gu, Li-Xin Zhang}
\date{}
\begin{document}
\maketitle

\begin{abstract}
This paper explores the Law of the Iterated Logarithm (LIL) for $m$-dependent sequences under the framework of sub-linear expectations. We first extend existing LIL results to sequences of independent, non-identically distributed random variables under sub-linear expectations. This extension serves as a crucial intermediary step, facilitating the subsequent establishment of the LIL for $m$-dependent stationary sequences. On the other hand, we also establish necessary conditions for $m$-dependent sequences in sub-linear expectation spaces. \end{abstract}

\section{Introduction}
The study of probability limit theorems under sub-linear expectations \cite{peng08, peng19} has garnered significant attention due to its ability to model uncertainty beyond classical linear frameworks. Sub-linear expectations generalize traditional probabilistic expectations by accommodating model uncertainty and allowing for the representation of a range of possible probability measures. This framework is particularly useful in areas such as finance \cite{peng07}, insurance \cite{Gianin06}, and robust statistics \cite{peng20}, where uncertainty and variability are inherent.

One of the fundamental results in probability theory is the Law of the Iterated Logarithm (LIL), which describes the fluctuations of partial sums of random variables \cite{HW41, wittman}. Extending the LIL to sub-linear expectation spaces involves overcoming challenges related to the non-additive nature of capacities. Recent advancements have established the LIL for independent random variables under sub-linear expectations \cite{chen14}, even when these variables are not identically distributed \cite{zhang21a}. These results provide a deeper understanding of the asymptotic behavior of sums of random variables in uncertain environments.

In parallel, the concept of $m$-dependence plays a crucial role in the analysis of dependent random variables. An $m$-dependent sequence allows each random variable to depend on only a fixed number of preceding variables, making it a versatile model for various applications, including time series analysis \cite{sen65}, spatial statistics \cite{Park09}, and network theory \cite{graph14}. The LIL for $m$-dependent sequences has been extensively studied in classical probability settings, where it has been proven under various dependence structures and moment conditions \cite{Petrov1984, chen97}. These studies demonstrate the robustness of the LIL in accommodating dependencies among random variables, thereby extending its applicability beyond independent sequences.

Proving the LIL in the context of sub-linear expectations and $m$-dependent sequences presents unique challenges. One of the primary difficulties arises from the discontinuity of capacities, which are inherently non-additive and may lack certain regularity properties present in linear expectations. This discontinuity complicates the analysis of partial sums, especially when dealing with truncated variables. Additionally, approximating the square of the sum of truncated variables is non-trivial, which needs the second order Choquet integral condition $C_{\V}(X^2)<\infty$. 

In this paper, we extend the LIL in \cite{zhang22note} to sequences of independent random variables that are not identically distributed under sub-linear expectations. This generalization serves as a crucial intermediary, facilitating the subsequent proof of the LIL for $m$-dependent stationary sequences by leveraging the independence structure within $m$-dependent frameworks. Additionally, we derive necessary conditions for $m$-dependent stationary sequences, suggesting that the assumptions in our LIL result could potentially be weakened.

The paper is structured as follows: Section \ref{notations} introduces key notations and preliminary concepts in sub-linear expectations, setting up the framework for our analysis. Section \ref{main_result} presents our main results, including the LIL for both independent and $m$-dependent sequences. Section \ref{inequalities} presents the technical tools and properties needed to support the main theorems, while Section \ref{proofs} provides detailed proofs, addressing challenges such as capacity discontinuity and the approximation of truncated sums.

\input{notation}
\input{main_result}
\input{lemma}

\input{proofs}
\input{references.bbl}

\end{document}

%% file: notation.tex
\section{ Settings and Notations.}\label{notations}
\setcounter{equation}{0}

We use the framework and notations of Peng \cite{peng08, peng19}. Let  $(\Omega,\mathcal F)$
 be a given measurable space  and let $\mathcal{H}$ be a linear space of real measurable functions
defined on $(\Omega,\mathcal F)$ such that if $X_1,\ldots, X_n \in \mathcal{H}$  then $\varphi(X_1,\ldots,X_n)\in \mathcal{H}$ for each
$\varphi\in C_{l,Lip}(\mathbb R^n)$,  where $C_{l,Lip}(\mathbb R^n)$ denotes the linear space of (local Lipschitz)
functions $\varphi$ satisfying
\begin{eqnarray*} & |\varphi(\bm x) - \varphi(\bm y)| \le  C(1 + |\bm x|^m + |\bm y|^m)|\bm x- \bm y|, \;\; \forall \bm x, \bm y \in \mathbb R^n,&\\
& \text {for some }  C > 0, m \in \mathbb  N \text{ depending on } \varphi. &
\end{eqnarray*}
$\mathcal{H}$ is considered as a space of ``random variables''.   We also denote $C_{b,Lip}(\mathbb R^n)$ the space of bounded  Lipschitz
functions.

\begin{define}\label{def1.1} A  sub-linear expectation $\mbE$ on $\mathcal{H}$  is a function $\mbE: \mathcal{H}\to \overline{\mathbb R}$ satisfying the following properties: for all $X, Y \in \mathcal H$, we have
\begin{description}
  \item[\rm (a)]  Monotonicity: If $X \ge  Y$ then $\mbE [X]\ge \mbE [Y]$;
\item[\rm (b)] Constant preserving: $\mbE [c] = c$;
\item[\rm (c)] Sub-additivity: $\mbE[X+Y]\le \mbE [X] +\mbE [Y ]$ whenever $\mbE [X] +\mbE [Y ]$ is not of the form $+\infty-\infty$ or $-\infty+\infty$;
\item[\rm (d)] Positive homogeneity: $\mbE [\lambda X] = \lambda \mbE  [X]$, $\lambda\ge 0$.
 \end{description}
 Here $\overline{\mathbb R}=[-\infty, \infty]$, $0\cdot \infty$ is defined to be $0$. The triple $(\Omega, \mathcal{H}, \mbE)$ is called a sub-linear expectation space. Given a sub-linear expectation $\mbE $, let us denote the conjugate expectation $\mbe$of $\mbE$ by
$$ \mbe[X]:=-\mbE[-X], \;\; \forall X\in \mathcal{H}. $$
\end{define}
 By Theorem 1.2.1 of Peng \cite{peng19}, there exists a family of finite additive linear expectations $E_{\theta}: \mathcal{H}\to \overline{\mathbb R}$ indexed by $\theta\in \Theta$, such that
\begin{equation}\label{linearexpression} 
\mbE[X]=\max_{\theta\in \Theta} E_{\theta}[X] \; \text{ for } X \in \mathcal{H} \text{ with } \mbE[X] \text{ being finite}. 
\end{equation}
Moreover, for each $X\in \mathcal{H}$, there exists $\theta_X\in \Theta$ such that $\mbE[X]=E_{\theta_X}[X]$ if $\mbE[X]$ is finite.

\begin{define}

\begin{description}
  \item[ \rm (i)] ({\em Identical distribution}) Let $\bm X_1$ and $\bm X_2$ be two $n$-dimensional random vectors defined
respectively in sub-linear expectation spaces $(\Omega_1, \mathcal{H}_1, \mbE_1)$
  and $(\Omega_2, \mathcal{H}_2, \mbE_2)$. They are called identically distributed, denoted by $\bm X_1\overset{d}= \bm X_2$  if
$$ \mbE_1[\varphi(\bm X_1)]=\mbE_2[\varphi(\bm X_2)], \;\; \forall \varphi\in C_{b,Lip}(\mathbb R^n). $$
 A sequence $\{X_n;n\ge 1\}$ of random variables is said to be identically distributed if $X_i\overset{d}= X_1$ for each $i\ge 1$.
\item[\rm (ii)] ({\em Independence})   In a sub-linear expectation space  $(\Omega, \mathcal{H}, \mbE)$, a random vector $\bm Y =
(Y_1, \ldots, Y_n)$, $Y_i \in \mathcal{H}$ is said to be independent to another random vector $\bm X =
(X_1, \ldots, X_m)$ , $X_i \in \mathcal{H}$ under $\mbE$  if for each test function $\varphi\in C_{b,Lip}(\mathbb R^m \times \mathbb R^n)$
we have
$ \mbE [\varphi(\bm X, \bm Y )] = \mbE \big[\mbE[\varphi(\bm x, \bm Y )]\big|_{\bm x=\bm X}\big].$
 \item[\rm (iii)] ({\em Independent random variables}) A sequence of random variables $\{X_n; n\ge 1\}$
 is said to be independent, if  $X_{i+1}$ is independent to $(X_{1},\ldots, X_i)$ for each $i\ge 1$.
 \item[\rm (iv)](m-dependence) A sequence of random variables(or random vectors) $\{X_n;n\geq 1\}$ is said to be $m$-dependent if there exists an integer $m$ such that for every $n$ and every $j\geq m+1$, $(X_{n+m+1},\cdots,X_{n+j})$ is independent of $(X_1,\cdots,X_n)$. In particular, if $m=0$, $\{X_n;n\geq 1\}$ is an independent sequence.
\item[\rm (v)](stationary) A sequence of random variables(or random vectors) $\{X_n;n\geq 1\}$ is said to be stationary if for every positive integers $n$ and $p$, $(X_1,\cdots,X_n)\overset{d}{=}(X_{1+p},\cdots,X_{n+p})$.
\item[\rm (vi)](linear stationary) A sequence of random variables(or random vectors) $\{X_n;n\geq 1\}$ is said to be linear stationary if for every positive integers $n$ and $p$, $X_1+\cdots+X_n\overset{d}{=}X_{1+p}+\cdots+X_{n+p}$.
 \end{description}
\end{define}

 Let $(\Omega, \mathcal{H}, \mbE)$ be a sub-linear expectation space.  We denote   $(\V,\mv)$ to be a pair of  capacities with the properties that
\begin{equation}\label{eq1.3}
   \mbE[f]\le \V(A)\le \mbE[g]\;\;
\text{ if } f\le I_A\le g, f,g \in \mathcal{H} \text{ and } A\in \mathcal F,
\end{equation}
 $$ \V \text{ is sub-additive in sense that } \V(A\bigcup B)\le \V(A)+\V(B)  \text{ for all } A,B\in \mathcal F  $$
and $\mv(A):= 1-\V(A^c)$, $A\in \mathcal F$.
We call $\V$ and $\mv$ the upper and the lower capacity, respectively. In general, we can choose $\V$ as
\begin{equation}\label{eq1.5} \V(A):=\inf\{\mbE[\xi]: I_A\le \xi, \xi\in\mathcal{H}\},\;\; \forall A\in \mathcal F.
\end{equation}
To distinguish this   capacity from others, we denote it by $\hV$, and $\mhV(A)=1-\hV(A)$. $\hV$ is the largest capacity satisfying \eqref{eq1.3}.

When there exists  a family of probability measure on $(\Omega,\mathcal{F})$ such that
\begin{equation}\label{eq1.7} \mbE[X]=\sup_{P\in \mathcal{P}}P[X]=:\sup_{P\in \mathcal{P}}\int  XdP ,
\end{equation} $\V$ can be defined as
\begin{equation}\label{eq1.6} \V(A)=\sup_{P\in \mathcal{P}}P(A).
\end{equation}
We denote this capacity by $\V^{\mathcal{P}}$, and $\mv^{\mathcal{P}}(A)=1-\V^{\mathcal{P}}(A)$.

 Also, we define the  Choquet integrals/expecations $(C_{\V},C_{\mv})$  by
$$ C_V[X]=\int_0^{\infty} V(X\ge t)dt +\int_{-\infty}^0\left[V(X\ge t)-1\right]dt $$
with $V$ being replaced by $\V$ and $\mv$ respectively.
If $\V_1$ on the sub-linear expectation space $(\Omega_1,\mathcal{H}_1,\mbE_1)$ and $\V_2$ on the sub-linear expectation space $(\Omega_2,\mathcal{H}_2,\mbE_2)$  are two capacities having the property \eqref{eq1.3}, then for any random variables $X_1\in \mathcal{H}_1$ and $\tilde X_2\in \mathcal{H}_2$ with $X_1\overset{d}=\tilde X_2$, we have
\begin{equation}\label{eqV-V}
\V_1(X_1\ge x+\epsilon)\le \V_2(\tilde X_2\ge x), \;\;  \mv_1(X_1\ge x+\epsilon)\le \mv_2(\tilde X_2\ge x) \text{ for all } \epsilon>0 \text{ and } x
\end{equation}
and
$C_{\V_1}(X_1)=C_{\V_2}(\tilde X_2).
$

\begin{define}\label{def3.1}
   A function $V:\mathcal F\to [0,1]$ is called to be  countably sub-additive if
$$ V\Big(\bigcup_{n=1}^{\infty} A_n\Big)\le \sum_{n=1}^{\infty}V(A_n) \;\; \forall A_n\in \mathcal F. $$
\end{define}

Since the capacity $\hV$ defined as in \eqref{eq1.5}  may be not countably sub-additive, we  consider its countably sub-additive extension.
\begin{define}   A  countably sub-additive extension $\Vstar$  of $\hV$   is defined by
\begin{equation}\label{outcapc} \Vstar(A)=\inf\Big\{\sum_{n=1}^{\infty}\hV(A_n): A\subset \bigcup_{n=1}^{\infty}A_n\Big\},\;\; \mvstar(A)=1-\Vstar(A^c),\;\;\; A\in\mathcal F,
\end{equation}
where $\hV$ is defined as in \eqref{eq1.5}.
\end{define}

As shown in Zhang \cite{zhang16}, $\Vstar$ is countably sub-additive, $\Vstar(A)\le \hV(A)$ and,  if $V$ is also a sub-additive (resp. countably sub-additive) capacity satisfying
\begin{equation}\label{eqVprop} V(A)\le \mbE[g]  \text{ whenever }  I_A\le g\in \mathcal{H},
\end{equation}
 then $V(A)\le \hV$ (resp. $V(A)\le \Vstar(A)$.

 \begin{define} Another countably sub-additive capacity generated by $\mbE$ can be defined as follows:
 \begin{equation}\label{tildecapc} \tildeCapc(A)=\inf\Big\{\lim_{n\to\infty}\mbE[\sum_{i=1}^n g_i]: I_A\le \sum_{n=1}^{\infty}g_n, 0\le g_n\in\mathcal{H}\Big\},\;\;\; A\in\mathcal F.
\end{equation}
\end{define}

 $\tildeCapc$ is  countably sub-additive and  has the property \eqref{eqVprop}, and so, $\tildeCapc(A)\le \Vstar(A)$. Furthermore,
 if $\mbE$ has the form \eqref{eq1.7}, then
$$ \V^{ \mathcal{P}}(A)=\sup_{P\in \mathcal{P}}P(A)\le \tildeCapc(A)\le \Vstar(A), \;\; A\in\mathcal{F}.  $$

Finally, we introduce the condition (CC) proposed in Zhang \cite{zhang21a}. We say that the sub-linear expectation $\mbE$ satisfies the condition (CC) if
\begin{equation}
    \mbE[X]=\sup_{P\in\mathcal{P}}P[X],\enspace X\in\mathcal{H}_b,
    \label{CC}
\end{equation}
	where $\mathcal{H}_b=\{f\in\mathcal{H}:f\enspace is\enspace bounded\}$, $\mathcal{P}$ is a countable-dimensionally weakly compact family of probability measures on $(\Omega,\sigma(\mathcal{H}))$ in sense that, for any $Y_1, Y_2,\cdots\in\mathcal{H}_b$ and any sequence $\{P_n\}\subset\mathcal{P}$, there is a subsequence $\{n_k\}$ and a probability measure $P\in\mathcal{P}$ for which
$$
	\lim_{k\rightarrow\infty}P_{n_k}[\varphi(Y_1,\cdots,Y_d)]=P[\varphi(Y_1,\cdots,Y_d)],\enspace \varphi\in C_{b,Lip}(\mathbb{R}^d),d\geq1.
$$

We denote
\begin{equation}
	\V^{\mathcal{P}}(A)=\sup_{P\in\mathcal{P}}P(A),\enspace A\in\sigma(\mathcal{H}),
\end{equation}
and it is obvious that $\V^{\mathcal{P}}\leq\widehat{\V}^*\leq \hV$. Let
\begin{align*}
 \mathcal{P}^e=& \left\{P: P \text{ is a probability measure on }\sigma(\mathcal{H}) \right. \\
 & \left.\text{ such that }
P[X]\le \mbE[X] \text{ for all } X\in \mathcal{H}_b \right\}.
\end{align*}

Zhang \cite{zhang21b} has shown the following three statements are equivalent: (i) the condition (CC) is satisfied with some $\mathcal{P}$; (ii) the condition (CC) is satisfied with  $\mathcal{P}^e$; (iii) $\mbE$ is regular in the sense that $\mbE[X_n]\searrow 0$ whenever $\mathcal{H}_b\ni X_n\searrow 0$. It shall be mentioned that  $\sup_{P\in\mathcal{P}_1}P[X] =\sup_{P\in\mathcal{P}_2}P[X]\; \forall  X\in\mathcal{H}_b$ does not imply $\sup_{P\in\mathcal{P}_1}P(A) =\sup_{P\in\mathcal{P}_2}P(A)\; \forall  A\in\sigma(\mathcal{H})$.

Through this paper, for real numbers $x$ and $y$, denote $x\vee y=\max(x,y)$, $x\wedge y=\min(x,y)$, $x^+=\max(0,x)$, $x^-=\max(0,-x)$ and $\log x=\ln \max(e,x)$. For a random variable $X$, because $XI\{|X|\le c\}$  may be not in $\mathcal{H}$, we will truncate it in the form $(-c)\vee X\wedge c$ denoted by $X^{(c)}$. We  denote
$$\breve{\mathbb E}[X]=\lim_{c\to \infty} \mbE[X^{(c)}]$$
if the limit exists.

%% file: main_result.tex
\section{Main Results}
\label{main_result}
We first present a law of the iterated logarithm for independent, though not identically distributed, sequences within the framework of sub-linear expectation probability space $\sles$. The first theorem serves as a mutual generalization of Theorem 2.2 in Zhang \cite{zhang22note}.

\begin{theorem}
    Let $\{X_n;n=1,2,\cdots\}$ be a sequence of independent random variables under $\sles$. Denote $s_n^2=\sum_{i=1}^n\mbE[X_i^2]$, $\underline{s}_n^2=\sum_{i=1}^n\mbe[X_i^2]$, $t_n=\sqrt{2\log\log s_n^2}$ and $a_n=s_nt_n$. Suppose that condition (CC) is satisfied. Assume that 
    \begin{equation}
        \frac{\underline{s}_n^2}{s_n^2}\to r > 0,
    \label{ratio}
    \end{equation}
    \begin{equation}
        s_n^2\to\infty, \quad |X_n|\leq \alpha_n\frac{s_n}{t_n} \text{ with } \alpha_n\to 0
        \label{bound}
    \end{equation}
    and
    \begin{equation}
        \frac{\sum_{j=1}^n|\mbE[X_j]|+|\mbe[X_j]|}{a_n}\to0.
        \label{mean-condition}
    \end{equation}
     are satisfied, then for $V = \V^{\mathcal{P}}, \mathbb{C}^*$ or $\Vstar$,
    \begin{equation}
        V\left(C\left\{\frac{\sum_{i=1}^n X_i}{\sqrt{2s_n^2\log\log s_n^2}}\right\}=[-\sqrt{t},\sqrt{t}]\right)=
        \begin{cases}
            1 & \text{ if } t\in[r, 1],\\
            0 & \text{ if } t\notin[r,1].
        \end{cases}
    \end{equation}
    \label{th-upper}
\end{theorem}

Now, we turn our attention to the $m$-dependent case. The following theorem presents the law of the iterated logarithm for $m$-dependent and stationary random variables.
\begin{theorem}
Let $\{X_n;n=1,2,\cdots\}$ be a m-dependent stationary sequence in a sub-linear space $\sles$ satisfying following conditions :
\begin{itemize}
    \item[(1)] $\bE[X_1] = \be[X_1] = 0$;
    \item[(2)] $\bE[X_1^2]<\infty$;
    \item[(3)] $C_{\V}\left[\frac{X_1^2}{\log\log|X_1|}\right]<\infty$.  
\end{itemize}
Then we have 
\begin{equation}
\mvstar\left(\lsigma\le\limsup_{n\to\infty}\frac{S_n}{\sqrt{2n\log\log n}}\le\usigma\right)=1\label{LIL}
\end{equation}
for some $0\le \lsigma\le\usigma<\infty$ where the following limits exist:
\begin{equation}
    \usigma^2:=\lim_{n\to\infty}\frac{\bE[S_n^2]}{n},\quad \lsigma^2:=\lim_{n\to\infty}\frac{\be[S_n^2]}{n}.
    \label{limit-sigma}
\end{equation}
\label{Th-lower}
\end{theorem}
\begin{remark}
    Due to the countably subadditivity of $\Vstar$, the conditions in Theorem \ref{Th-lower} are equivalent to those in the i.i.d case(c.f. Theorem 5.5 in Zhang \cite{zhang21a}).
\end{remark}
Next theorem shows the exact lower bound for $m$-dependent sequences.
\begin{theorem}
\label{th:m-dependt}
    Let $\{X_n;n=1,2,\cdots\}$ be a sequence of $m$-dependent and linear stationary random variables in a sub-linear space $\sles$ with 
    \begin{equation}
        C_{\V}(X_1^2)<\infty,
    \end{equation}
    and 
    \begin{equation}
        \bE[X_1] = \be[X_1] = 0,
        \label{m-mean-condition}
    \end{equation}
    Denote $S_n=\sum_{i=1}^nX_i$ and $\lsigma^2, \usigma^2$ are defined as in (\ref{limit-sigma}). Further assume $\lsigma^2>0$ and $\mbE$ satisfies the condition (CC), then for $V = \V^{\mathcal{P}}, \mathbb{C}^*$ or $\Vstar$,
    \begin{equation}
        V\left(C\left\{\frac{\sum_{i=1}^n X_i}{\sqrt{2n\log\log n}}\right\}=[-\sigma,\sigma]\right)=
        \begin{cases}
            1 & \text{ if } \sigma\in[\lsigma, \usigma],\\
            0 & \text{ if } \sigma\notin[\lsigma, \usigma].
        \end{cases}
        \label{upperLIL}
    \end{equation}
\end{theorem}

Finally, the following theorem establishes the necessary condition for the law of the iterated logarithm.

\begin{theorem}
\label{th:necessary}
    Let $\{X_n;n=1,2,\cdots\}$ be a sequence of $m$-dependent and linear stationary random variables in a sub-linear space $\sles$ satisfying condition (CC). If there exists $M$ such that 
    \begin{equation}
        \Vstar\left(\limsup_{n\to\infty}\frac{|S_n|}{\sqrt{2n\log\log n}}\geq M \right) < 1,
        \label{eq: unbounded}
    \end{equation}
    then 
    \begin{equation}
        C_{\V}\left[\frac{X_1^2}{\log\log |X_1|}\right]<\infty
        \label{choquet-bound}
    \end{equation}
    and
    \begin{equation}\label{eq: necessary mean}
        \lim_{n\to\infty}\frac{\bE[S_n]}{n}=\lim_{n\to\infty}\frac{\be[S_n]}{n}=0.
    \end{equation}
\end{theorem}

%% file: lemma.tex
\section{Related inequalities and properties}\label{inequalities}
The first lemma establishes the properties of random variables for which  $C_{\V}\left[\frac{X^2}{\log\log|X|}\right]<\infty$. The specific details and proof of these properties can be found in Lemma 6.1 of Zhang \cite{zhang21a}.
\begin{lemma} \label{lem:property} Suppose $X\in \mathcal{H}$.
\begin{description}
  \item[\rm (i)]
 For any $\delta>0$,
$$ \sum_{n=1}^{\infty} \V\big(|X|\ge \delta \sqrt{n\log\log n} \big)<\infty \;\; \Longleftrightarrow C_{\V}\left[\frac{X^2}{\log\log|X|}\right]<\infty.
$$
 \item[\rm (ii)]
   If $C_{\V}\left[\frac{X^2}{\log\log|X|}\right]<\infty$, then for any $\delta>0$ and $p>2$,
$$ \sum_{n=1}^{\infty} \frac{\mbE\big[\big(|X|\wedge (\delta \sqrt{n\log\log n})\big)^p\big]}{(n\log\log n)^{p/2}}<\infty. $$
 \item[\rm (iii)] $C_{\V}\left[\frac{X^2}{\log\log|X|}\right]<\infty$, then for any $\delta>0$,
 $$ \mbE[X^2\wedge (2\delta n\log\log n)]=o(\log \log n) $$
 and $$ \breve{\mathbb E}[(|X|-\delta \sqrt{2  n\log\log n})^+]=o(\sqrt{\log\log n/n}).\;\;  $$
\end{description}
\end{lemma}

Next lemma gives the exponential inequalities and the proof is in Lemma 3.1 of \cite{zhang22note}.
\begin{lemma}\label{lem:ExpIneq}   Suppose that  $\{X_1,\ldots, X_n\}$ is a sequence  of independent random variables
on $\sles$. Set
$A_n(p,y)=\sum_{i=1}^n\mbE[(X_i^+\wedge y )^p]$ and
$ \breve{B}_{n,y}=\sum_{i=1}^n \mbE[(X_i\wedge y)^2]$.   Then  for all $p\ge 2$, $x,y>0$, $0<\delta\le 1$,
\begin{align}
    & \hV\Big( \max_{k\le n}  \sum_{i=1}^k(X_i-\bE[X_i])\ge x\Big)
    \le    \hV\big(\max_{k\le n} X_k> y \big)
    +\exp\left\{-\frac{x^2}{2(xy+\breve{B}_{n,y})   }\right\}.
    \label{eq:ExpIneq.1}
\end{align}
and
\begin{align}\label{eq:ExpIneq.2}
& \hV\Big( \max_{k\le n}  \sum_{i=1}^k(X_i-\bE[X_i])\ge x\Big) \nonumber\\
\le &  \hV\big(\max_{k\le n} X_k> y \big)
  +2\exp\{p^p\}\Big\{\frac{A_n(p,y)}{y^p} \Big\}^{\frac{\delta x}{10y}}
+\exp\left\{-\frac{x^2}{2\breve{B}_{n,y}(1+\delta)   }\right\}.
\end{align}
\end{lemma}

The following self-normalized law of the iterated logarithm for martingales  is  Lemma 13.8 of   Pe\~{n}a et al. \cite{pena09}.
\begin{lemma}\label{lem:selfLILM} Let $\{X_n; n\ge 1\}$ be a martingale difference sequence in a probability space $(\Omega,\mathcal{F},P)$ with respect to the
filtration $\{\mathcal F_n\}$ such that $|X_n|\le  m_n$ a.s. for some $\mathcal F_{n-1}$-measurable random variable
$m_n$, with $U_n\to \infty$  and $m_n/\{U_n(\log\log U_n)^{-1/2}\}\to 0$  a.s., where $U_n^2=\sum_{i=1}^n X_i^2$.  Then
$$ \limsup_{n\to \infty}\frac{\sum_{i=1}^n X_i}{U_n(2\log\log U_n)^{1/2}}=1\;\; a.s. $$
\end{lemma}

The following lemma gives the relation between the sub-linear expectation  and the conditional expectation under a probability, the  proof of which  can be found in  Guo and Li \cite{Guo&Li21} (see also Hu et al. \cite{Hu21} and Guo et al. \cite{GLL22}.

\begin{lemma}\label{lem3} Let $\{X_n;n\ge 1\}$ be a sequence of independent random variables in the sub-linear expectation space $\sles$ with \eqref{CC}. Denote
$$ \mathcal{F}_n=\sigma(X_1,\ldots,X_n) \; \text{ and }\; \mathcal{F}_0=\{\emptyset,\Omega). $$
Then for each $P\in\mathcal{P}$, we have
$$ E_P\left[\varphi(X_n)|\mathcal{F}_{n-1}\right]\le \mbE[\varphi(X_n)]\;\; a.s., \;\; \varphi\in C_{b,Lip}(\mathbb R). $$
\end{lemma}

The following lemma is a special case of Lemma3.1 in Guo\cite{Guo22} for identically distributed random variables with $C_{\V}[X^2_1]<\infty$. 
\begin{lemma}
	Let $\{X_n;n\geq 1\}$ be a sequence of identically distributed random variables in the sub-linear expectation space $\sles$ with $C_{\V}[X_1^2]<\infty$. Then there exists a sequence of positive constants $\{\kappa_i;i\leq1\}$ with $\kappa_i\downarrow0$ and $\kappa_i\sqrt{i/\log\log i}\rightarrow\infty$ such that
\begin{equation*}
	\sum_{i=16}^\infty\frac{\bE\left[\left(\vert X_i\vert-\kappa_i\sqrt{i/\log\log i}\right)^+\right]}{\sqrt{2i\log\log i}}<\infty.
\end{equation*}
\label{truncation}
\end{lemma}

The next lemma is the Borel-Cantelli lemma for a countably sub-additive capacity,  the proof of which can be referred to \cite{zhang21a}.
\begin{lemma}\label{lem:BCdirect} Let $V$ be a   countably sub-additive capacity and $\sum_{n=1}^{\infty}V(A_n)<\infty$. Then
$$ V(A_n\; i.o.)=0, \;\; \text{ where } \{ A_n\; i.o.\}=\bigcap_{n=1}^{\infty}\bigcup_{i=n}^{\infty}A_i. $$
\end{lemma}

The next lemma is  the converse part of the Borel-Cantelli lemma (see \cite{zhang21b} and \cite{gu23}).
\begin{lemma}\label{lem:BCconverse} Let $\sles$ be a sub-linear expectation  space with a capacity $\V$ having the property (\ref{eq1.3}). Suppose the condition (CC) is satisfied.
Suppose that $\{X_n;n\ge 1\}$ is a sequence of independent random variables in  $\sles$.
 \begin{description}
   \item[\rm (i)] If $\sum_{n=1}^{\infty}\mv(X_n>1)<\infty$,  then
$$\mv^{\mathcal{P}}\left( X_n>1 \; i.o.\right)=0.
$$
   \item[\rm (ii)]   If $\sum_{n=1}^{\infty}\V(X_n\ge 1)=\infty$, then
$$ \V^{\mathcal{P}}\left( X_n\ge 1\;\; i.o. \right)=1.
$$
 \end{description}
\end{lemma}

The last two lemmas give the relationship between $\mbE$ and $\bE$.
The first about $\bE$ is Proposition 1.1 of Zhang \cite{zhang21b}.
\begin{lemma}\label{lemma3.6} Consider a subspace of $\mathcal{H}$ as
\begin{equation}\label{eqlem3.6.1} \mathcal{H}_1=\big\{X\in\mathcal{H}: \lim_{c,d\to \infty}\mbE\big[(|X|\wedge d-c)^+\big]=0\big\}.
\end{equation}
Then for any $X\in\mathcal{H}_1$, $\bE[X]$ is well defined, and $(\Omega,\mathcal{H}_1,\bE)$ is a sub-linear expectation space.
\end{lemma}
\begin{lemma} 
Suppose $X,Y\in \mathcal{H}_1$ and $X$ is independent to $Y$ under the sub-liner expectation $\mbE$. If $\bE[X]=\be[X]=0$ and $\bE[|X|]<\infty, \bE[|Y|]<\infty$, then 
\begin{equation}
    \bE[XY]=0.
\end{equation}
\label{l1}
\end{lemma}
\begin{proof}
    For any $b>0$, note 
    \begin{equation*}
        |(XY)^{(b)}-X^{(\sqrt{b})}Y^{(\sqrt{b})}|\leq(|X|\wedge \sqrt{b})\cdot(|Y|-\sqrt{b})^++(|Y|\wedge \sqrt{b})\cdot(|X|-\sqrt{b})^+. 
    \end{equation*}
    It follows from the independence and $X,Y\in \mathcal{H}_1$ that
    \begin{align*}
       & \left|\mbE\left[(XY)^{(b)}\right]-\mbE\left[X^{(\sqrt{b})}Y^{(\sqrt{b})}\right]\right|\\
        \leq& \mbE\left[(|X|\wedge \sqrt{b})\cdot(|Y|-\sqrt{b})^+\right] + \mbE\left[(|Y|\wedge \sqrt{b})\cdot(|X|-\sqrt{b})^+\right]\\
        = & \mbE\left[|X|\wedge\sqrt{b}\right]\cdot\mbE\left[(|Y|-\sqrt{b})^+\right] + \mbE\left[|Y|\wedge\sqrt{b}\right]\cdot\mbE\left[(|X|-\sqrt{b})^+\right]\to 0
    \end{align*}
    as $b\to\infty$. By the independence again, we have
    \begin{align*}
        \left|\bE\left[XY\right]\right| &= \lim_{b\to\infty}\left|\mbE\left[(XY)^{(b)}\right]\right|\\
        & =\lim_{b\to\infty}\left|\mbE\left[X^{(\sqrt{b})}Y^{(\sqrt{b})}\right]\right|\\
        & \le \lim_{b\to\infty}\mbE\left[|Y|\wedge \sqrt{b}\right]\cdot\left(\left|\mbE\left[X^{(\sqrt{b})}\right]\right|+\left|\mbe\left[X^{(\sqrt{b})}\right]\right|\right)\\
        & = \bE[|Y|]\cdot\left(|\bE[X]|+|\be[X]|\right)= 0 .
    \end{align*}
\end{proof}

%% file: proofs.tex
\section{Proofs of the theorems}\label{proofs}
In this section, we will prove the theorems presented in Section \ref{main_result}. To demonstrate Theorem \ref{th-upper}, we will utilize the self-normalized law of the iterated logarithm for martingale difference sequences (Lemma \ref{lem:selfLILM}) and first establish the following three propositions. Combining Proposition \ref{prop2} and \ref{prop3}, we can immediately obtain Theorem \ref{th-upper}.

\begin{proposition}
Let $\{X_n;n=1,2,\cdots\}$ be a sequence of independent random variables in a sub-linear space $\sles$ satisfying $|X_n|\leq\alpha_n\frac{s_n}{t_n}$ for some $\alpha_n\to0$ and $s_n^2\to\infty$. Denote $\underline{s}_n^2=\sum_{i=1}^n\mbe[X_i^2]$ and $V_n^2=\sum_{i=1}^nX_i^2$, and assume that there exists a real number $r$ such that
\begin{equation}
    \frac{\underline{s}_n^2}{s_n^2}\to r.
\end{equation}
Then 
\begin{equation}
    \mvstar\left(r\leq\liminf_{n\to\infty}{\frac{V_n^2}{s_n^2}}\leq\limsup_{n\to\infty}{\frac{V_n^2}{s_n^2}}\leq1\right)=1.
\end{equation}
\label{prop1}
\end{proposition}
\begin{proof}
    By Lemma 3.3 of Wittmann \cite{wittman}, for $\forall\lambda>1$, there exists a subsequence $\{n_k\}\subset\mathbb{N}$ with 
    $$\lambda a_{n_k}\leq a_{n_{k+1}}\leq \lambda^3 a_{n_k+1}.$$
    It can be checked that 
    $$\log s_{n_{k+1}^2}\sim\frac12\log a_{n_{k+1}}\geq ck.$$
    By Lemma \ref{lem:ExpIneq}, denote $b_n=\alpha_n\frac{s_n}{t_n}$, let $x = \varepsilon s_{n_k}^2$ and $y=2b_{n_k}^2$, we have
    \begin{align*}
        &\hV\left(\max_{n\leq n_k}\sum_{i=1}^n(X_i^2-\mbE[X_i^2])\geq\varepsilon s_{n_k}^2\right)\\
        &\leq \exp\left\{-\frac{\varepsilon^2 s_{n_k}^4}{2(2\varepsilon s_{n_k}^2b_{n_k}^2+b_{n_k}^2\sum_{i=1}^{n_k}\mbE[X_i^2])}\right\}\\
        &\leq\exp\left\{-\frac{\varepsilon^2t_{n_k}^2}{2(2\varepsilon + 1)\alpha_{n_k}^2}\right\}\\
        &\leq\exp\{-2\log\log s_{n_k}^2\}\leq ck^{-2}.
    \end{align*}
    Hence
    $$\sum_{k=1}^\infty \hV\left(\max_{n\leq n_k}\sum_{i=1}^n(X_i^2-\mbE[X_i^2])\geq\varepsilon s_{n_k}^2\right)<\infty.$$
    It follows that 
    $$\mvstar\left(\limsup_{k\to\infty}\frac{\max_{n\leq n_k}\sum_{i=1}^n(X_i^2-\mbE[X_i^2])}{s_{n_k}^2}\leq 0\right)=1$$
    by the countable sub-additivity of $\Vstar$, $\Vstar(A)\leq\hV(A)$ and the Borel-Cantelli lemma-Lemma \ref{lem:BCdirect}.
    It follows from the above equation that
    $$\mvstar\left(\limsup_{n\to\infty}\frac{\sum_{i=1}^nX_i^2}{s_n^2}\leq 1\right).$$
    On the other hand, replacing $X_i^2$ with $-X_i^2$ yields
    $$\mvstar\left(\limsup_{k\to\infty}\frac{\max_{n\leq n_k}\sum_{i=1}^n(X_i^2-\mbe[X_i^2])}{s_{n_k}^2}\geq 0\right)=1,$$
    which means
    $$\mvstar\left(\limsup_{n\to\infty}\frac{\sum_{i=1}^nX_i^2}{s_n^2}\geq r\right)$$
    since $\frac{\underline{s}_n^2}{s_n^2}\to r$.
\end{proof}

\begin{proposition}
    Under the conditions of Proposition \ref{prop1},  we further assume that the condition (CC) is satisfied. Then for any $t\in[r,1]$, we have 
    \begin{equation}
        \V^{\mathcal{P}}\left(\lim_{n\to\infty}\frac{V_n^2}{s_n^2}=t\right)=1.
    \end{equation}
    \label{prop2}
\end{proposition}
\begin{proof}
    By (\ref{linearexpression}), there are finite additive linear expectations $E_{j,1}$ and $E_{j,2}$ such that
    $$E_{j,1}[X_j^2]=\mbe[X_j^2] \quad\text{and}\quad E_{j,2}[X_j^2]=\mbE[X_j^2].$$
    For any $t\in [r,1]$, we can write $t=\alpha r + (1-\alpha)$. Let $E_j=\alpha E_{j,1} + (1-\alpha) E_{j,2}$, $E_j$ is a finite additive linear expectation with $E_j\leq \mbE$. 

    By Proposition 2.1 of \cite{zhang21b}, we can find a new sub-linear space $(\tilde{\Omega}, \tilde{\mathcal{H}},\tilde{\mathbb{E}})$ defined on a metric space $\tilde{\Omega}=\mathbb{R}^\infty$, with a copy $\{\tilde{X}_n;n\geq 1\}$ on $(\tilde{\Omega}, \tilde{\mathcal{H}},\tilde{\mathbb{E}})$ of $\{X_n;n\geq 1\}$ and a probability measure $Q$ on $\tilde{\Omega}$ such that $\{\tilde{X}_n;n\geq 1\}$ is a sequence of independent random variables under $Q$, satisfying
    \begin{equation}
        E_Q\left[\varphi(\tilde{X}_i)\right] = E_i\left[\varphi(\tilde{X}_i)\right]\quad \text{for all } \varphi\in C_{b,Lip}(\mathbb{R}),
        \label{eqQ}
    \end{equation}
    \begin{equation*}
    E_Q\left[\varphi(\tilde{X}_1,\cdots,\tilde{X}_d)\right]\leq \tilde{\mathbb{E}}\left[\varphi(\tilde{X}_1,\cdots,\tilde{X}_d)\right]=\mbE\left[\varphi(X_1,\cdots,X_d)\right]\quad \text{for all } \varphi\in C_{b,Lip}(\mathbb{R}^d)\nonumber
    \end{equation*}
    and 
    \begin{equation}
        \tilde{v}(B)\leq Q(B)\leq\tilde{V}(B)\quad \text{for all } B\in \sigma(\tilde{Y_1},\tilde{Y_2},\cdots).
        \label{eqv}
    \end{equation}
    where
    \begin{equation*}
        \tilde{V}(A)=\tilde{\mathbb{V}}^{\tilde{\mathcal{P}}}(A)=\sup_{P\in\tilde{\mathcal{P}}}P(A) \text{ and } \tilde{v}(A) = 1-\tilde{V}(A^c), A\in \tilde{\mathcal{F}},
    \end{equation*}
    and $\tilde{\mathcal{P}}$ is the family of all probability measures $P$ on $(\tilde{\Omega},\tilde{\mathcal{F}})$ with the property
    \begin{equation*}
        E_P[\varphi]\leq \tilde{\mathbb{E}}[\varphi] \text{ for bounded } \varphi\in\tilde{\mathcal{H}}.
    \end{equation*}
    Notice $|X_i|\leq b_i$ and by (\ref{eqQ}), it follows that
    \begin{equation*}
        E_Q[\tilde{X}_j^4]=E_j[X_j^4]\leq b_j^2E_j[X_j^2]\leq b_j^2\mbE[X_j^2]. 
    \end{equation*}
    Apply exponential inequality again, it follows that
    \begin{align*}
        &Q\left(\max_{n_k+1\leq n\leq n_{k+1}}\left|\sum_{j=n_k+1}^n(\tilde{X}_j^2-E_Q[\tilde{X}_j^2])\right|\geq\varepsilon s_{n_{k+1}}^2\right)\\
        &\leq 2\exp\left\{-\frac{\varepsilon^2 s_{n_k}^4}{2(2\varepsilon s_{n_k}^2b_{n_k}^2+b_{n_k}^2\sum_{i=1}^{n_k}\mbE[X_i^2])}\right\}\\
        &\leq2\exp\left\{-\frac{\varepsilon^2t_{n_k}^2}{2(2\varepsilon + 1)\alpha_{n_k}^2}\right\}\\
        &\leq2\exp\{-2\log\log s_{n_k}^2\}\leq ck^{-2}.
    \end{align*}
    Hence for all $\varepsilon>0$,
    \begin{equation*}
        \sum_{k=1}^{\infty}Q\left(\max_{n_k+1\leq n\leq n_{k+1}}\left|\sum_{j=n_k+1}^n(\tilde{X}_j^2-E_Q[\tilde{X}_j^2])\right|\geq\varepsilon s_{n_{k+1}}^2\right)<\infty.
    \end{equation*}
    Then there exists a sequence $0<\varepsilon \searrow 0$ such that 
    \begin{equation*}
        \sum_{k=1}^{\infty}Q\left(\max_{n_k+1\leq n\leq n_{k+1}}\left|\sum_{j=n_k+1}^n(\tilde{X}_j^2-E_Q[\tilde{X}_j^2])\right|\geq\varepsilon_k s_{n_{k+1}}^2\right)<\infty.
    \end{equation*}
    By (\ref{eqv}), we have
    \begin{equation*}
        \sum_{k=1}^{\infty}\tilde{v}\left(\max_{n_k+1\leq n\leq n_{k+1}}\left|\sum_{j=n_k+1}^n(\tilde{X}_j^2-E_Q[\tilde{X}_j^2])\right|\geq\varepsilon_k s_{n_{k+1}}^2\right)<\infty.
    \end{equation*}
    Since $(Y_1,\cdots, Y_{n_{k+1}})\overset{d}{=}(\tilde{Y}_1,\cdots,\tilde{Y}_{n_{k+1}})$, by (\ref{eq1.6}) we have
    \begin{equation*}
        \sum_{k=1}^{\infty}\mathcal{V}^{\mathcal{P}}\left(\max_{n_k+1\leq n\leq n_{k+1}}\left|\sum_{j=n_k+1}^n(X_j^2-E_Q[\tilde{X}_j^2])\right|\geq2\varepsilon_k s_{n_{k+1}}^2\right)<\infty.
    \end{equation*}
    By Borel-Cantelli lemma-Lemma \ref{lem:BCconverse}, it follows that
    \begin{equation*}
        \mathcal{V}^{\mathcal{P}}(A_k\enspace i.o.)=0\text{ with }A_k=\left\{\frac{\max_{n_k+1\leq n\leq n_{k+1}}\left|\sum_{j=n_k+1}^n(X_j^2-E_Q[\tilde{X}_j^2])\right|}{s_{n_{k+1}}^2}>2\varepsilon_k\right\}.
    \end{equation*}
    On the event $(A_k\enspace i.o.)^c$,
    \begin{equation*}
        \lim_{k\to\infty}\frac{\max_{n_k+1\leq n\leq n_{k+1}}\left|\sum_{j=n_k+1}^n(X_j^2-E_Q[\tilde{X}_j^2])\right|}{s_{n_{k+1}}^2}=0,
    \end{equation*}
    which implies
    \begin{align*}
        &\lim_{n\to\infty}\frac{\sum_{i=1}^n X_i^2}{s_n^2}=\lim_{n\to\infty}\frac{\sum_{i=1}^n E_Q[\tilde{X}_i^2]}{s_n^2}=\lim_{n\to\infty}\frac{\sum_{i=1}^n E_i[X_i^2]}{s_n^2}\\
        &=\lim_{n\to\infty}\frac{\alpha\sum_{i=1}^n \mbe[X_i^2]+(1-\alpha)\sum_{i=1}^n \mbE[X_i^2]}{s_n^2}=\alpha r+(1-\alpha)=t.
    \end{align*}
    Hence $\V^{\mathcal{P}}\left(\lim_{n\to\infty}\frac{V_n^2}{s_n^2}=t\right)=1$.
\end{proof}

\begin{proposition}
    Let $\{X_n;n=1,2,\cdots\}$ be a sequence of independent random variables in a sub-linear space $\sles$. Suppose there is a family of probability measures $\mathcal{P}$ on $(\Omega, \sigma(\mathcal{H}))$ such that the sub-linear expectation $\mbE$ satisfies condition (CC). Assume that (\ref{ratio}) - (\ref{mean-condition}) are satisfied.
    Then
    \begin{equation}
        \mv^{\mathcal{P}}\left(C\left\{\frac{\sum_{i=1}^n X_i}{V_n\sqrt{2\log\log s_n^2}}\right\}=[-1,1]\right)=1.
        \label{cluster}
    \end{equation}
    \label{prop3}
\end{proposition}

\begin{proof}
    Let $P\in\mathcal{P}, \mathcal{F}_n=\sigma(X_1, \cdots,X_n),\mathcal{F}_0=\{\emptyset, \Omega\}$. By Lemma \ref{lem3}, it follows that
    \begin{equation*}
        E_P[X_i|\mathcal{F}_{i-1}]\leq\mbE[X_i]\text{ and } E_P[-X_i|\mathcal{F}_{i-1}]\leq\mbE[-X_i],
    \end{equation*}
    which implies
    \begin{equation}
        \mbe[X_i]\leq E_P[X_i|\mathcal{F}_{i-1}]\leq\mbE[X_i].
        \label{martingale-mean}
    \end{equation}
    Denote $U_n^2:=\sum_{i=1}^n(X_i-E_P[X_i|\mathcal{F}_{i-1}])^2$ and $b_n=\alpha_n\sqrt{\frac{s_n}{t_n}}$.
    By Proposition \ref{prop1}, we have $V_n^2\sim s_n^2 \enspace a.s. P$. It follows from (\ref{mean-condition}) and (\ref{martingale-mean}) that
    \begin{align*}
        |V_n^2-U_n^2|&=\left|\sum_{i=1}^nE_P[X_i|\mathcal{F}_{i-1}](2X_i-E_P[X_i|\mathcal{F}_{i-1}])\right|\\
        &\leq 3b_n\sum_{i=1}^n|E_P[X_i|\mathcal{F}_{i-1}]|\\
        &\leq 3b_n\sum_{i=1}^n(|\mbE[X_j]|+|\mbe[X_j]|)=o(s_n^2).
    \end{align*}
    Hence $U_n^2\sim s_n^2\enspace a.s.P$ and $b_n^2/(U_n(\log\log U_n)^{-1/2})=\alpha_n\to0$. By Lemma \ref{lem:selfLILM}, we have
    \begin{align*}
        &P\left(\limsup_{n\to\infty}\frac{\sum_{i=1}^n(X_i-E_P[X_i|\mathcal{F}_{i-1}])}{V_n\sqrt{2\log\log s_n^2}}=1\right)\\
        =& P\left(\limsup_{n\to\infty}\frac{\sum_{i=1}^n(X_i-E_P[X_i|\mathcal{F}_{i-1}])}{U_n\sqrt{2\log\log U_n}}=1\right)=1.
    \end{align*}
     From (\ref{mean-condition}) and (\ref{martingale-mean}) again we yield
     \begin{equation*}
         \frac{\sum_{i=1}^nE_P[X_i|\mathcal{F}_{i-1}]}{V_n\sqrt{2\log\log s_n^2}}\to0\enspace a.s.P,
     \end{equation*}
     and then 
     \begin{equation*}
         P\left(\limsup_{n\to\infty}\frac{\sum_{i=1}^nX_i}{V_n\sqrt{2\log\log s_n^2}}=1\right)=1.
     \end{equation*}
     Replacing $X_i$ with $-X_i$, we obtain (\ref{cluster}).
\end{proof}

To prove Theorem \ref{Th-lower}, we utilize the method outlined in \cite{chen97} along with the LIL for the i.i.d. case.
\begin{proof}[\textbf{Proof of Theorem \ref{Th-lower}}]
    First, we show that the limit of $\bE[S_n^2]/n$ exists. By Lemma \ref{l1}, we have
\begin{align*}
    \bE[S_n^2]\leq n\bE[X_1^2] + 2\sum_{i=1}^n\sum_{j=i+1}^{n\wedge(i+m)}\bE[X_iX_j]\leq nC.
\end{align*}
\par With the same argument of CLT for $m$-dependent sequence under sub-linear expectations \cite{gu23CLT}, it follows that the limits of $\bE[S_n^2]/n$ and $\be[S_n^2]/n$ exist. We denote 
$$
 \lsigma^2:=\lim_{n\to\infty}\frac{\be[S_n^2]}{n},\usigma^2:=\lim_{n\to\infty}\frac{\bE[S_n^2]}{n} 
$$
and we have $0\leq\lsigma^2\leq\usigma^2<\infty$.
\par Now, we assume that the sequence $\{X_n\}$ is 1-dependent under $\mbE$. Fix some integer $p>1$, let
$$
Y_k = \sum_{(k-1)p<j<kp}X_j,\quad k\ge 1,\quad T_n = \sum_{k=1}^n Y_k.
$$
\par Since $\{X_n\}$ is 1-dependent and stationary, both $\{Y_k;k=1,2,\cdots\}$ and $\{X_{kp};k=1,2,\cdots\}$ are i.i.d. sequences under $\mbE$. It is easy to verify that 
$$
0=\sum_{j=1}^{kp-1}\be[X_j]\le \be[Y_1]\le\bE[Y_1]\le\sum_{j=1}^{kp-1}\bE[X_j]=0
$$
and
$$
C_{\V}\left[\frac{Y_1^2}{\log\log|Y_1|}\right]<\infty.
$$
\par By Theorem 5.5 of \cite{zhang21a}, it follows that 
\begin{align}
   & \mvstar\left(\lsigma_p\leq\limsup_{k\to \infty}\frac{T_n}{\sqrt{2n\log\log n}}\leq \usigma_p\right)=1, \label{eqY}\\
   & \mvstar\left(\lsigma_X\leq\limsup_{k\to \infty}\frac{\sum_{k=1}^n X_{kp}}{\sqrt{2n\log\log n}}\leq \usigma_X\right)=1,\label{eqX}
\end{align}
where we denote $\lsigma_p^2:=\be[Y_1^2] =\be[S_{p-1}^2]$, $\usigma_p^2:=\bE[Y_1^2]=\bE[S_{p-1}^2]$, $\lsigma_X^2:=\be[X_1^2]$ and $\usigma_X^2:=\bE[X_1^2]$. Thus, it follows from \eqref{eqY} and \eqref{eqX} that 
\begin{align*}
    \mvstar\left(\lsigma_p-\lsigma_X\leq\limsup_{k\to \infty}\frac{\sum_{k=1}^{np} X_{k}}{\sqrt{2n\log\log n}}\leq \usigma_p+\usigma_X\right)=1.
\end{align*}
\par Note 
$$
\mvstar\left(\lim_{n\to\infty}\frac{|X_n|}{\sqrt{2n\log\log n}}=0\right)=1,
$$
it is obvious that
$$
    \left|\limsup_{n\to\infty}\frac{S_n}{\sqrt{2n\log\log n}}-\limsup_{n\to\infty}\frac{S_{np}}{\sqrt{2np\log\log (np)}}\right|\leq \limsup_{n\to\infty}\frac{|S_p|}{\sqrt{2n\log\log n}}=0
$$
under $\mvstar$. Thus, for any integer $p$, 
\begin{align*}
    \mvstar\left(p^{-1/2}(\lsigma_p-\lsigma_X)\leq\limsup_{k\to \infty}\frac{\sum_{k=1}^{n} X_{k}}{\sqrt{2n\log\log n}}\leq p^{-1/2}(\usigma_p+\usigma_X)\right)=1.
\end{align*}
Note that 
\begin{align*}
    p^{-1/2}\usigma\to 0,\quad
    \frac{\lsigma_p^2}{p} = \frac{\mbe[S_{p-1}^2]}{p}\to \lsigma^2,\quad
    \frac{\usigma_p^2}{p} = \frac{\mbE[S_{p-1}^2]}{p}\to \usigma^2,
\end{align*}
and by the countably subadditivity of $\Vstar$, we have
\begin{align*}
    &\Vstar\left(\limsup_{k\to\infty}\frac{\sum_{k=1}^{n} X_{k}}{\sqrt{2n\log\log n}}<\lsigma\text{ or } \limsup_{k\to\infty}\frac{\sum_{k=1}^{n} X_{k}}{\sqrt{2n\log\log n}}> \usigma\right)\\
    =&\Vstar\left(\bigcup_{p=1}^\infty\left\{\limsup_{k\to\infty}\frac{\sum_{k=1}^{n} X_{k}}{\sqrt{2n\log\log n}}<\frac{\lsigma_p-\lsigma_X}{p^{1/2}}\text{ or } \limsup_{k\to\infty}\frac{\sum_{k=1}^{n} X_{k}}{\sqrt{2n\log\log n}}> \frac{\usigma_p+\usigma_X}{p^{1/2}}\right\}\right)\\
    \leq & \sum_{p=1}^\infty \Vstar\left(\limsup_{k\to\infty}\frac{\sum_{k=1}^{n} X_{k}}{\sqrt{2n\log\log n}}<\frac{\lsigma_p-\lsigma_X}{p^{1/2}}\text{ or } \limsup_{k\to\infty}\frac{\sum_{k=1}^{n} X_{k}}{\sqrt{2n\log\log n}}> \frac{\usigma_p+\usigma_X}{p^{1/2}}\right)=0.
\end{align*}

Thus we obtain the equation \eqref{LIL} for 1-dependent sequence. Now we consider m-dependent case. Let $\{X_n\}$ be a m-dependent and stationary sequence, denote 
$$
Z_k=\sum_{j=(k-1)m+1}^{km}X_j,
$$
then $\{Z_k\}$ is 1-dependent. It is easy to verify that the conditions hold for $\{Z_k\}$, thus, we have
\begin{align*}
     &\mvstar\left(\lsigma\le\limsup_{n\to\infty}\frac{S_n}{\sqrt{2n\log\log n}}\le\usigma\right)\\
  =&\mvstar\left(m^{1/2}\lsigma\le\limsup_{n\to\infty}\frac{\sum_{k=1}^nZ_k}{\sqrt{2n\log\log n}}\le m^{1/2}\usigma\right)=1.
\end{align*}
   
\end{proof}

\begin{proof}[\textbf{Proof of Theorem \ref{th:m-dependt}}]
Since $\{X_n;n=1,2,\cdots\}$ is stationary and $C_{\hV}(|X_1|^2)<\infty$, by Lemma \ref{truncation}, there exists a sequence of positive numbers $\alpha_n$ such that $\alpha_n\searrow0$, $\alpha_n\sqrt{n/\log\log n}\to\infty$ and
    \begin{equation*}
        \sum_{i=16}^n\frac{\bE\left[\left(|X_i|-\alpha_i\sqrt{i/\log\log i}\right)^+\right]}{\sqrt{2i\log\log i}}<\infty.
    \end{equation*}
    There exists a sequence $M_i\nearrow\infty$ such that 
    \begin{equation*}
        \sum_{i=16}^n\frac{M_i\bE\left[\left(|X_i|-\alpha_i\sqrt{i/\log\log i}\right)^+\right]}{\sqrt{2i\log\log i}}<\infty.
    \end{equation*}
    Denote $b_i=\alpha_i\sqrt{i/\log\log i}$ and $Y_i=X_i^{(b_i)}$, it follows from Kronecker Lemma that for any $P\in\mathcal{P}$,
    \begin{equation}
        \lim_{n\to\infty}\frac{\sum_{i=16}^nM_i\bE\left[\left(|X_i|-\alpha_i\sqrt{i/\log\log i}\right)^+\right]}{\sqrt{2n\log\log n}}=0\text{ and }\lim_{n\to\infty}\frac{\sum_{i=1}^n|X_i-Y_i|}{\sqrt{2n\log\log n}}=0 \enspace a.s.P.
        \label{M_i}
    \end{equation}
    Define $h_0=0, h_n=h_{n-1}+l_n+m$, where $l_n$ is chosen such that $l_n\nearrow\infty$
    \begin{equation}
        l_n\leq M_{h_{n-1}+1}^{1/4} \quad\text{ and } \quad l_n\alpha_n\to 0.
        \label{l_n-construction}
    \end{equation}
    Denote 
    \begin{equation*}
        Z_n = \sum_{i=h_{n-1}+1}^{h_n-m}Y_i, \quad W_n=\sum_{i=h_n-m+1}^{a_n}X_i.
    \end{equation*}
    Then $\{Z_n;n\geq1\}$ and $\{W_n;n\geq1\}$ are sequences of independent random variables, respectively. Furthermore, $\{W_n;n\geq1\}$ is an identical distributed sequence. Now we need to show the former part converges to our target while the latter part vanishes. 

    First, we need to estimate $s_n^2=\sum_{k=1}^n\mbE[Z_n^2]$. On one hand, we have
    \begin{align*}
        &\mbE[Z_k^2]-\bE\left[\left(\sum_{i=h_{k-1}+1}^{h_k-m}X_i\right)^2\right]\\
        \leq& \bE\left[\left(\sum_{i=h_{k-1}+1}^{h_k-m}Y_i\right)^2-\left(\sum_{i=h_{k-1}+1}^{h_k-m}X_i\right)^2\right]\\
        \leq &\sum_{i=h_{k-1}+1}^{h_k-m} \bE\left[Y_i^2-X_i^2\right] + \sum_{1\leq|i-j|\leq m}\bE[Y_iY_j-X_iX_j]+\sum_{|i-j|>m}\bE[Y_iY_j-X_iX_j]\\
        \leq & \sum_{i=h_{k-1}+1}^{h_k-m} \bE\left[|Y_i^2-X_i^2|\right] + \sum_{1\leq|i-j|\leq m}\bE[|Y_iY_j-X_iX_j|]+\sum_{|i-j|>m}\mbE[Y_iY_j]=: I_k + II_k+ III_k,
    \end{align*}
    where the last inequality is implied from Lemma \ref{l1}. For $I_k$, we have
    \begin{align}
        \limsup_{n\to\infty}\frac{\sum_{k=1}^nI_k}{h_n}\leq\limsup_{n\to\infty}\frac{\sum_{i=1}^{h_n}\bE[(X_1^2-b_i^2)^+]}{h_n}\leq C_{\V}[(X_1^2-C)^+]\to0
        \label{I_k}
    \end{align}
    as $C\to\infty$. For $II_k$, we have
    \begin{align*}
        II_k=&\sum_{1\leq|i-j|\leq m}\bE[|Y_iY_j-X_iX_j|]\\
        \leq& \sum_{1\leq|i-j|\leq m}\left(\bE[|Y_iY_j-X_iY_j|] + \bE[|X_iY_j-X_iX_j|]\right)\\
        \leq & \sum_{1\leq|i-j|\leq m}\left(\bE[Y_j^2]^{1/2}\bE[(X_i-Y_i)^2]^{1/2}+\bE[X_i^2]^{1/2}\bE[(X_j-Y_j)^2]^{1/2}\right)\\
        \leq& C_m\sum_{i=h_{k-1}+1}^{h_k-m}\bE[(X_i-Y_i)^2]^{1/2}
        \leq C_m\sum_{i=h_{k-1}+1}^{h_k-m}\bE[(X_1^2-b_i^2)^+]^{1/2},
    \end{align*}
    where the second inequality is derived from the Cauchy-Schwarz inequality and the third inequality is due to $\bE[X_1^2]<\infty$. Hence 
    \begin{equation}
        \limsup_{n\to\infty}\frac{\sum_{k=1}^nII_k}{h_n}\leq\limsup_{n\to\infty}C_m\frac{\sum_{i=1}^{h_n}\bE[(X_1^2-b_i^2)^+]^{1/2}}{h_n}\leq  C_mC_{\V}[(X_1^2-C)^+]^{1/2}\to0
    \end{equation}
    as $C\to\infty$. Finally, for $III_k$, note that for $j-i>m$, $Y_j$ is independent to $Y_i$ and $\mbE[Y_iY_j]\leq\mbE\left[|Y_i|\right]\left(\left|\mbE[Y_j]\right|+\left|\mbe[Y_j]\right|\right)$, it follows that
    \begin{align*}
        III_k=&\sum_{|i-j|>m}\mbE[Y_iY_j]\leq 2\sum_{i=h_{k-1}+1}^{h_k-m}\sum_{j=i+m+1}^{h_k-m}\mbE[|Y_i|]\left(\left|\mbE[Y_j]\right|+\left|\mbe[Y_j]\right|\right)\\
        \leq& \sum_{i=h_{k-1}+1}^{h_k-m}\mbE[|Y_i|]\sum_{j=h_{k-1}+1}^{h_k-m}\left(\left|\mbE[Y_j]\right|+\left|\mbe[Y_j]\right|\right)\\
        \leq & l_n\bE[|X_1|]\sum_{j=h_{k-1}+1}^{h_k-m}\left(\left|\mbE[Y_j]-\bE[X_j]\right|+\left|\mbe[Y_j]-\be[X_j]\right|\right)\\
        \leq & Cl_n\sum_{i=h_{k-1}+1}^{h_k-m}\bE\left[|X_i-Y_i|\right]\\
        \leq & C\sum_{i=h_{k-1}+1}^{h_k-m} M_i\bE\left[(|X_1|-b_i)^+\right],
    \end{align*}
    where the third inequality is derived from (\ref{m-mean-condition}) and the last inequality is obtained by the construction of $l_n$ in (\ref{l_n-construction}). Then it follows from (\ref{M_i}) that
    \begin{equation}
        \limsup_{n\to\infty}\frac{\sum_{k=1}^nIII_k}{h_n}\leq \limsup_{n\to\infty}C\frac{\sum_{i=1}^{h_n} M_i\bE\left[(|X_1|-b_i)^+\right]}{h_n}=0.
        \label{III_k}
    \end{equation}
    Combining (\ref{I_k})-(\ref{III_k}), we can conclude 
    \begin{equation*}
        \limsup_{n\to\infty}\frac{\sum_{k=1}^n\mbE[Z_k^2]-\sum_{k=1}^n\bE\left[\left(\sum_{i=h_{k-1}+1}^{h_k-m}X_i\right)^2\right]}{h_n}\leq0.
    \end{equation*}
    On the other hand, we can similarly show 
    \begin{equation*}
        \limsup_{n\to\infty}\frac{\sum_{k=1}^n\bE\left[\left(\sum_{i=h_{k-1}+1}^{h_k-m}X_i\right)^2\right]-\sum_{k=1}^n\mbE[Z_k^2]}{h_n}\leq0.
    \end{equation*}
    Note 
    \begin{align*}
        &\lim_{n\to\infty}\frac{\sum_{k=1}^n\bE\left[\left(\sum_{i=h_{k-1}+1}^{h_k-m}X_i\right)^2\right]}{h_n}=\lim_{n\to\infty}\frac{\bE\left[\left(\sum_{i=h_{n-1}+1}^{h_n-m}X_i\right)^2\right]}{l_n+m}\\
        =&\lim_{n\to\infty}\frac{\bE\left[\left(\sum_{i=1}^{l_n}X_i\right)^2\right]}{l_n}\cdot\frac{l_n}{l_n+m}=\usigma^2.
    \end{align*}
    Hence we have
    \begin{equation}
        \lim_{n\to\infty}\frac{s_n^2}{h_n}=\lim_{n\to\infty}\frac{\sum_{k=1}^n\mbE[Z_k^2]}{h_n}=\lim_{n\to\infty}\frac{\sum_{k=1}^n\bE\left[\left(\sum_{i=h_{k-1}+1}^{h_k-m}X_i\right)^2\right]}{h_n}=\usigma^2.
    \end{equation}
    Similarly, we can verify that $\underline{s}_n^2=\sum_{i=1}^n\mbe[Z_i^2]\sim h_n\lsigma^2$ and $\frac{\underline{s}_n^2}{s_n^2}\to\frac{\lsigma^2}{\usigma^2}=:r$.
    The expression for $Z_n$ is bounded as follows:
    \begin{equation*}
        |Z_n|\leq\sum_{i=h_{n-1}+1}^{h_n-m}b_i\leq l_n\alpha_{h_n}\sqrt{\frac{h_n}{\log\log h_n}}\leq l_n\alpha_{n}\sqrt{\frac{h_n}{\log\log h_n}}.
    \end{equation*}
    Additionally, it is obvious that $s_n^2\to\infty$ since $\usigma^2\ge\lsigma^2>0$ and it follows from (\ref{l_n-construction}) that $l_n\alpha_n\to0$. Consequently, (\ref{bound}) is satisfied.

    Moreover, denote $d_n=\sqrt{2n\log\log n}$ and $a_n=\sqrt{2s_n^2\log\log s_n^2}\sim\usigma d_{h_n}$, it follows from (\ref{m-mean-condition}) and (\ref{M_i}) that
    \begin{align}
        \frac{\sum_{i=1}^n\left|\mbE[Z_i]\right| + \left|\mbe[Z_i]\right|}{a_n}=&\frac{\sum_{i=1}^n\left|\mbE[Z_i]-\bE\left[\sum_{j={h_{i-1}+1}}^{h_i-m}X_j\right]\right| + \left|\mbe[Z_i]-\be\left[\sum_{j={h_{i-1}+1}}^{h_i-m}X_j\right]\right|}{a_n}\nonumber\\
        \leq & \frac{\sum_{i=1}^n\sum_{j={h_{i-1}+1}}^{h_i-m}\bE\left[(|X_j|-b_j)^+\right]}{a_n}\nonumber\\
        \leq &C\frac{\sum_{i=1}^{h_n}\bE\left[(|X_j|-b_j)^+\right]}{\sqrt{2h_n\log\log h_n}} \to 0.
        \label{eq:mean for z}
    \end{align}
    Applying Theorem \ref{th-upper}, we have for $V = \V^{\mathcal{P}}, \mathbb{C}^*$ or $\Vstar$,
    \begin{equation}
        V\left(C\left\{\frac{\sum_{i=1}^n Z_i}{d_{h_n}}\right\}=[-\sigma,\sigma]\right)=
        \begin{cases}
            1 & \text{ if } \sigma\in[\lsigma, \usigma],\\
            0 & \text{ if } \sigma\notin[\lsigma, \usigma].
        \end{cases}
        \label{eq:ZLIL-any}
    \end{equation}
    For $\{W_n;n\geq 1\}$, note $W_n$ is independent and identically distributed, we can apply Theorem 5.5 in \cite{zhang21a} and obtain 
    \begin{equation*}
        \mvstar\left(-\bE[W_1^2]^{1/2}\leq\liminf_{n\to\infty}\frac{\sum_{i=1}^nW_i}{\sqrt{2n\log\log n}}\leq\limsup_{n\to\infty}\frac{\sum_{i=1}^nW_i}{\sqrt{2n\log\log n}}\leq\bE[W_1^2]^{1/2}\right)=1.
    \end{equation*}
    Since $l_n\to\infty$ and $n=o(h_n)$, we have
    \begin{equation}
        \mvstar\left(\lim_{n\to\infty}\frac{\sum_{i=1}^nW_i}{d_{h_n}}=0\right)=1.
        \label{W-omit}
    \end{equation}
    For $h_n+1\leq N\leq h_{n+1}$, we define 
    \begin{equation*}
        T_N=\sum_{i=1}^n\left(\sum_{j=h_{i-1}+1}^{h_i-m}Y_j+\sum_{j=h_i-m+1}^{h_i}X_j\right)+\sum_{j=h_n+1}^N\left(Y_j1_{\{j\leq h_n-m\}}+X_j1_{\{j>h_n-m\}}\right).
    \end{equation*}
    Then from (\ref{W-omit}) we have
    \begin{equation}
        \Vstar\left(\lim_{n\to\infty}\left(\frac{T_{h_n}}{d_{h_n}}-\frac{\sum_{i=1}^nZ_i}{d_{h_n}}\right)\neq 0\right)=\Vstar\left(\lim_{n\to\infty}\frac{\sum_{i=1}^nW_i}{d_{h_n}}\neq0\right)=0.
        \label{W-omit2}
    \end{equation}
    Now by (\ref{M_i}), (\ref{eq:ZLIL-any}) and (\ref{W-omit2}), it is sufficient to show 
    \begin{equation}
        \Vstar\left(\lim_{n\to\infty}\max_{h_n+1\leq N\leq h_{n+1}}\left|\frac{T_N}{d_N}-\frac{T_{h_n}}{d_{h_n}}\right|\neq 0\right)=0.
        \label{target}
    \end{equation}
    For $h_n+1\leq N\leq h_{n+1}$, we have
    \begin{align*}
        \left|\frac{T_N}{d_N}-\frac{T_{h_n}}{d_{h_n}}\right|\leq\left|\frac{T_N-T_{h_n}}{d_N}+T_{h_n}\cdot\frac{d_{h_n}-d_N}{d_Nd_{h_n}}\right|\leq\frac{\sum_{j=h_n+1}^{h_{n+1}-m}|Y_j|}{d_{h_n}}+\frac{T_{h_n}}{d_{h_n}}\cdot \frac{d_{h_{n+1}}-d_{h_n}}{d_{h_n}}.
    \end{align*}
    For the first term, since $l_n\alpha_n\to 0$ and $h_{n+1}/h_n\to 1$, we have
    \begin{equation}
        \frac{\sum_{j=h_n+1}^{h_{n+1}-m}|Y_j|}{d_{h_n}}\leq\frac{l_{n+1}\alpha_{h_{n+1}}\sqrt{h_{n+1}/\log\log h_{n+1}}}{\sqrt{2h_{n}\log\log h_n}}\leq \frac{l_{n+1}\alpha_{n+1}}{\log\log h_n}\cdot\sqrt{\frac{h_{n+1}}{h_n}}\to0.
        \label{1st-term}
    \end{equation}
    For the second term, it follows that
    \begin{equation}
        \frac{d_{h_{n+1}}-d_{h_n}}{d_{h_n}}=\frac{d_{h_{n+1}}}{d_{h_n}}-1=\sqrt{\frac{h_{n+1}}{h_n}}\cdot \sqrt{\frac{\log\log h_{n+1}}{\log\log h_n}}-1\to0,
        \label{2nd-term}
    \end{equation}
    while $T_{h_n}/d_{h_n}$ is bounded for all $P\in\mathcal{P}$ by applying Kolmogorov's law of iterated logarithm(c.f. Corollary 5.2 in \cite{zhang21a}) to $\{Z_n;n\geq 1\}$. Combining (\ref{1st-term}) and (\ref{2nd-term}), we obtain (\ref{target}). 
    
    Finally, we show equation (\ref{upperLIL}). Since $\lsigma>0$ and $\frac{\underline{s}_n^2}{s_n^2}\to\frac{\lsigma^2}{\usigma^2}=r>0$, applying Theorem \ref{th-upper}, it follows that for $V = \V^{\mathcal{P}}, \mathbb{C}^*$ or $\Vstar$,
    \begin{equation}
        V\left(C\left\{\frac{\sum_{i=1}^n Z_i}{d_{h_n}}\right\}=[-\sigma,\sigma]\right)=
        \begin{cases}
            1 & \text{ if } \sigma\in[\lsigma, \usigma],\\
            0 & \text{ if } \sigma\notin[\lsigma, \usigma].
        \end{cases}
        \label{Z-main}
    \end{equation}
    Combining (\ref{W-omit2}), (\ref{target}) and (\ref{Z-main}), we obtain equation (\ref{upperLIL}). The proof is completed.
\end{proof}
\begin{proof}[\textbf{Proof of Theorem \ref{th:necessary}}]
    Suppose $C_{\V}\left[\frac{X_1^2}{\log\log |X_1|}\right]=\infty$, by Lemma \ref{lem:property}, we have for any $M >0$,
    \begin{equation*}
        \sum_{n=1}^{\infty}\hV(|X_1|\geq 2Md_n)=\infty.
    \end{equation*}
    There exists $M_n \nearrow \infty$ such that
    \begin{equation*}
        \sum_{n=1}^{\infty}\hV(|X_1|\geq 2M_nd_n)=\infty.
    \end{equation*}
    By (\ref{eqV-V}), it follows that
    \begin{equation*}
        \sum_{i=1}^{\infty}\hV(|X_{i(m+1)}|\geq M_id_i)=\infty.
    \end{equation*}
    By the $m$-dependence, $\{X_{i(m+1)};i\geq 1\}$ are independent under $\mbE$. Hence, by Borel-Cantelli Lemma \ref{lem:BCconverse}, we have
    \begin{equation*}
        \V^{\mathcal{P}}(|X_{i(m+1)}|\geq M_id_i\enspace i.o.)=1.
    \end{equation*}
    On the event $\{|X_{i(m+1)}|\geq M_id_i\enspace i.o.\}$, we have
    \begin{align*}
        \infty &= \limsup_{i\to\infty}\frac{|X_{i(m+1)}|}{\sqrt{2i\log\log i}} = \sqrt{m+1}\limsup_{i\to\infty}\frac{|X_{i(m+1)}|}{\sqrt{2i(m+1)\log\log (i(m+1))}}\\
        &\leq \sqrt{m+1}\limsup_{n\to\infty}\frac{|X_n|}{d_n}=\sqrt{m+1}\limsup_{n\to\infty}\frac{|S_n-S_{n-1}|}{d_n}\leq2\sqrt{m+1}\limsup_{n\to\infty}\frac{|S_n|}{d_n}.
    \end{align*}
    Hence 
    \begin{equation*}
        \Vstar(\limsup_{n\to\infty}\frac{|S_n|}{d_n}=\infty)=1,
    \end{equation*}
    which is contradictory to condition (\ref{eq: unbounded}). It follows that \eqref{choquet-bound} holds and then $\bE[X_1]$ and $\be[X_1]$ exist and are finite. Since $\be[X_1]\leq\frac{\be[S_n]}{n}\leq\frac{\bE[S_n]}{n}\leq \bE[X_1]$, with the same approach as \cite{gu23}, we can show that $\underline{\mu}:=\lim_{n\to\infty}\frac{\be[S_n]}{n}$ and $\bar{\mu}:=\lim_{n\to\infty}\frac{\bE[S_n]}{n}$ exist and are finite.
    
    Further, it follows from the strong law of large numbers(c.f. Theorem 3.3 in \cite{gu23}) that for any $\underline{\mu}\le\mu\le\bar{\mu}$, 
    \begin{equation}\label{eq:5.23}
        \Vstar\left(\lim_{n\to\infty}\frac{S_n}{n}= \mu\right)=1.
    \end{equation}
    On the other hand, by \eqref{eq: unbounded}, there exist $0<\tau<1$ such that 
    \begin{equation}\label{eq5.24}
        \mvstar\left(\limsup_{n\to\infty}\frac{|S_n|}{d_n}<M\right)>1-\tau>0.
    \end{equation}
    Combining \eqref{eq:5.23} and \eqref{eq5.24}, we have for any $\underline{\mu}\le\mu\le\bar{\mu}$,
    \begin{equation*}
        \Vstar\left(\lim_{n\to\infty}\frac{S_n}{n}= \mu \text{ and } \limsup_{n\to\infty}\frac{|S_n|}{d_n}<M\right)>1-\tau>0.
    \end{equation*}
    It implies that 
    \begin{equation*}
        -\lim_{n\to\infty}\frac{Md_n}{n}\leq\lim_{n\to\infty}\frac{\be[S_n]}{n}\leq\mu\leq\lim_{n\to\infty}\frac{\bE[S_n]}{n}\leq \lim_{n\to\infty}\frac{Md_n}{n} .
    \end{equation*}
    Hence, the equation \eqref{eq: necessary mean} holds.
\end{proof}